\documentclass[12pts]{amsart}
\theoremstyle{plain} \newtheorem{Thm}{ }[section]
\title{Rank two vector bundles with canonical determinant}
\author{Montserrat Teixidor i Bigas}
\address{Mathematics Department, Tufts University, Medford
MA02155, USA}

\begin{document}

\begin{abstract}
Denote by $B^k_{2,K}$ the locus of vector bundles of rank two that
have canonical determinant and at least $k$ sections. We show that
for a generic curve of genus $g$, $B^k_{2,K}$ is non-empty and has
a component of the expected dimension if $g$ is sufficiently
large.
\end{abstract}

\maketitle

\begin{section}{Introduction}
Denote by $U(r,d)$ the moduli space of
stable vector bundles of rank $r$ and degree $d$ on a fixed curve
$C$ of genus $g$. The Brill-Noether loci
 $B^k_{r,d}$ are defined as the
subsets of $U(r,d)$ consisting of vector bundles with $k$
 sections (note that we use
$k$ rather than $k-1$ in the notation). These can be given
 a locally determinantal
scheme structure. As such the expected dimension of $B^k_{r,d}$
for a generic curve is given by the Brill-Noether number $\rho
=r^2(g-1)+1-k(k-d+r(g-1))$ and the expected singular locus is
$B^{k+1}_{d,r}$. The case $r=1$ is classical and all the expected
results hold for the generic curve. For higher rank though,
several cases are known  where $\rho$ is positive and the locus is
empty (cf \cite{BGN}) or $\rho$ is negative and the locus is
non-empty (use \cite{Mefl}) or the singular locus is larger than
expected (see \cite{T1}). We shall add to the list of
abnormalities by showing that $B^k_{2,2g-2}$ is sometimes of
dimension larger than expected for the generic curve and therefore
reducible (see also \cite{BF}).

Denote now by $U(r,L)$ the moduli space of stable vector bundles
of rank $r$ and determinant $L$. Define similarly $B^k_{r,L}$  as
the subset of $U(r,L)$ consisting of vector bundles with $k$
sections . For $k>r$ (and $d\le 2g-2$), one expects that the
geometry of $B^k_{r,L}$ will very much depend on the position of
$L$ itself as a point of various $B^j_{1,d}$. The case $L={ O}$ is
not very interesting as a stable vector bundle of degree zero
cannot have sections. It is then natural to look at the special
case when $ L=K$ is the canonical bundle. When $r=2$, this locus
has a natural structure as a symmetric degeneracy locus (see
\cite{Mu}, \cite{H},\cite{MkBN},\cite{BF},
 \cite{MkFano} ).
Using this description, the expected dimension of  $B^k_{2,K}$ is
given by $\rho^k _K=3g-3-{k+1\choose 2}$. By its definition, this
is the minimum dimension of any component of these loci.

One can again ask the question (\cite{MkBN}) of whether the locus is
non-empty for $\rho^k_K>0$ and of dimension $\rho^k_K$ and empty for
$\rho^k_K<0$. It was conjectured in \cite{BF} that this would be the case.
Some evidence for this conjecture using curves of small genus
was provided in \cite{R}.
The main purpose of this paper is to partly prove these
statements.

We want to mention that the study of the loci $B^k_{2,K}$, apart
from its intrinsic interest, seem to contribute a lot to our
knowledge of the geometry of the curves. Mukai for example showed
in \cite{Mkg=11}  that if $C$ is a curve of genus eleven, then
$W^7_{2,K}$ is a $K3$ surface containing the curve. Something can
be said also for curves of higher genus.

Moreover, as $B^k_{2,K}\subset B^k_{2,2g-2}$, a study of
$B^k_{2,K}$ will sometimes provide some information on
$B^k_{2,2g-2}$ (see \ref{corollari}).

Our main result is the following

\bigskip

\begin{Thm}
\label{Teorema} {\bf Theorem} Let $C$ be a generic curve of genus
$g$. If $k=2k_1$ is even, then $B^k_{2,K}$ is non-empty when
$$g\ge k_1^2 \ \ if \ k_1> 2,\
g\ge 5 \ \  if  \ k_1=2,\ g\ge 3 \ \ if \  k_1=1.$$
 Moreover, under these conditions,
it has a component of the right dimension $\rho ^k_K$.

If $k=2k_1+1$ is odd, then $B^k_{2,K}$ is non-empty when $g\ge
(k_1)^2+k_1+1$. Moreover, under these conditions, it has a
component of the right dimension.
\end{Thm}

\begin{Thm}
\label{corollari} {\bf Corollary} If $C$ is a generic curve of
genus $g,\ (k_1)^2\le g<2(k_1)^2-k_1$,
 then $B^{2k_1}_{2,2g-2}$
is non-empty and has a component of dimension larger than expected.

If $C$ is a generic curve of genus $g,\ (k_1)^2+k_1+1\le
g<2(k_1)^2+k_1$, then $B^{2k_1+1}_{2,2g-2}$ is non-empty and has a
component of dimension larger than expected.
\end{Thm}
\bigskip
The author is a member of the EAGER group "Vector bundles on
Algebraic Curves"

\end{section}
\begin{section}{Proof of existence}

We start by recalling some definitions and results that we shall later need.

\begin{Thm}
\label{cce} {\bf Chain of elliptic curves} Let $C_i,i=1...g$ be
elliptic curves. Let $P_i,Q_i$
 be generic points on $C_i$. Take
$C$ to be the curve obtained by identifying $Q_i$ and $P_{i+1}$.
\end{Thm}

In order to give linear series on elliptic curves it will be useful to
have the following

\begin{Thm}
\label{slce}{\bf Remark} One can define a subspace of dimension
$k$ of sections of a fixed line bundle $L$ of degree $d$ by
specifying the $k$ distinct (minimum) desired vanishings of a
basis of the subspace at two different points $P,Q$ so that the
sum of the corresponding vanishings at $P,Q$ is $d-1$. In the case
when $L={ O}(aP+(d-a)Q)$, two of the vanishings could be chosen to
be $a,d-a$ (adding up to $d$ rather than $d-1$). These are the
only two vanishings that can add to $d$ if $P,Q$ are generic.
\end{Thm}

\begin{proof}
Assume given a line bundle of degree $d$ on an elliptic curve $C$
and two points $P,Q$ on $C$. If $L$ is general, there is no
section $s$ of $L$ such that the orders of vanishing of $s$ at $P$
and $Q$ add up to $d$. If such a section exists, then $L={
O}(aP+(d-a)Q)$ for some $a\le d$. Moreover, if $P,Q$ are generic
(i.e. they don't differ in a torsion point of order at most $d$)
such an $a$ is unique.

On the other hand, there is always a section (defined up to a
constant) whose orders of vanishing at $P,Q$ are at least given
numbers $a,b$ satisfying $a+b=d-1$. If $L$ is generic, the
vanishings of this section at $P,Q$ are precisely $a,b$
respectively. If $L={ O}(aP+(d-a)Q)$, then the vanishings are $a,
d-a=b+1$ respectively.

Hence, our statement is proved.

 To simplify our task,
we shall always list the vanishing $u_j$ at $P$ in increasing order and
the vanishing $v_j$ at $Q$ in decreasing order.

A similar argument can be used to define sections of a given rank two vector bundle.
Now every vanishing at $P,Q$ may appear twice.
\end{proof}

  On a reducible curve, the equivalent of a linear series is
 a limit linear series
introduced by Eisenbud and Harris (cf \cite{EH}) for the rank one case and
later generalised to arbitrary rank in \cite{T1}. For the convenience of the
reader, we reproduce the definition here:

\begin{Thm}
\label{sll} {\bf Limit linear series} of rank $r$ dimension $k$
and degree $d$ on a chain of $M$ (not necessarily elliptic) curves
consists of data I, II below for which data III,IV exist
satisfying conditions a-c.

I) For every component $C_i$, a vector bundle $E_i$ of rank $r$ and degree
$d_i$ and a $k$-dimensional space of sections of $E_i$.

II) For every node obtained by gluing $Q_i$ and $P_{i+1}$ an isomorphism of
 the
projectivisations of the fibers $(E_i)_{Q_i}$ and $(E_{i+1})_{P_{i+1}}$.

III) A positive integer $a$.

IV) For every node obtained by gluing $Q_i$ and $P_{i+1}$, basis
$s_{Q_i}^t, s_{P_{i+1}}^t,\ t=1...k$ of the vector spaces
$V_i$ and $V_{i+1}$.

Subject to the conditions

a) $\sum _{i=1}^M d_i-r(M-1)a=d$

b)The orders of vanishing at $Q_i,P_{i+1}$ of the sections of the chosen
basis satisfy $s_{Q_i}^t+ s_{P_{i+1}}^t\ge a$

c)Sections of the vector bundles $E_i(-aP_i),E_i(-aQ_i)$
 are completely determined
by their value at the nodes.
\end{Thm}

{\bf Example} When $r=1$ (cf \cite{EH}) the gluing in II is
irrelevant. One can take $d_i=d$ and then $a=d$ too. Conditions a)
and c) are automatically satisfied and one needs only check
condition b).

For instance, the limit of the canonical series on a curve as in \ref{cce}
is the unique linear series of rank one dimension $g$ and degree $2g-2$.
One can check that this limit is given as follows:
$$E_i={ O}((2g-2i)Q_i+(2i-2)P_i).$$
The space of sections is the unique subspace of $H^0(
O((2g-2i)Q_i+(2i-2)P_i))$ whose vanishing $u^i_j,v^i_j$ at $P_i,\
Q_i$ are given by,
$$u^i_e=i-3+e,e<i,\ u^i_k=i-2+e,e\ge i,\
v^i_e=2g-i-e,e\le i,\  v^i_e=2g-i-e-1,e>i$$ (use Remark
\ref{slce})

\begin{proof}
We now turn to the proof of \ref{Teorema} The main fact that we
use in the proof is the following: As $B^k_{2,K}$ is defined as
the pull-back of a symmetric degeneracy locus, its dimension at
every point is at least $\rho ^k_K$. Even better, given a family
of curves ${\bf C}\rightarrow B$ one can define similarly ${\bf
B}^k_{2,K}$ parametrising pairs consisting of a curve in the
family and a vector bundle of rank two with canonical determinant
and at least $k$ sections. Then, dim(${\bf B}^k_{2,K}$ )$\ge $
dim$(B)+\rho ^k_K$ at every point. If we can find a curve $C_0$ in
such a family for which dim($B^k_{2,K})=\rho ^k_K$, then the
generic fiber of the map ${\bf B}^k_{2,K}\rightarrow B$ is
non-empty and has dimension $\rho ^k_K$. This is what we shall do
with a curve $C_0$ as in \ref{cce}.

We first
prove existence in the case of even $k=2k_1$ and under the assumption
$g\ge (k_1)^2$($g>(k_1)^2$ if $k_1\le 2$).

We must give a family of limit linear series on a curve as in \ref{cce} of rank two
dimension $k$ and canonical determinant.
We need to show that the dimension of this family is $3g-3-{k+1\choose 2}$ and that this
family is not part of a larger family of linear series with canonical determinant.

Consider the first $(k_1)^2$ components $C_i$. Write
$$i=a^2+2c+\epsilon ,\ 0\le a\le k_1-1,\ 0\le c\le a-1, \ \epsilon =1,2,\
\  or\ c=a,\ \epsilon =1$$ Take the vector bundle to be of the
form
$${ O}((g-i-c+a)Q+(c-a+i-1)P)\oplus { O}((g-i-a+c)Q+(a-c+i-1)P)$$
Note that the determinant of this vector bundle is ${
O}((2g-2i)Q+2(i-1)P)$ corresponding to the canonical bundle. Take
the space of sections whose vanishing $u^i_j,\ v^i_j$ at $P_i$,
$Q_i$ are given respectively by
$$u^i_{2e-1}=u^i_{2e}=i+e-a-3, \  v^i_{2e-1}=v^i_{2e}=g-i-e+a+1, e\le c$$
$$ u^i_{2c+1}=i+c-a-1 \ \if\ \epsilon =1, \
   u^i_{2c+1}=i+c-a-2 \ \ if\ \epsilon =2$$
$$u^i_{2c+2}=i+c-a-1 $$
$$v^i_{2c+1}=g-i-c+a,\ v^i_{2c+2}=g-i-c+a-1 \ \ if\ \epsilon =1, \
                       v^i_{2c+2}=g-i-c+a \ \ if\ \epsilon =2$$
$$u^i_{2e-1}=u^i_{2e}=i+e-a-2, \  v^i_{2e-1}=v^i_{2e}=g-i-e+a, c+1<e\le a$$
$$ i\not= a^2+2a+1,\ u^i_{2a+1}=i-c+a-1 \ \ if\ \epsilon =1, \
   u^i_{2a+1}=i-c+a-2 \ \ if\ \epsilon =2$$
$$u^i_{2a+2}=i+a-c-1 $$
$$v^i_{2a+1}=g-i+c-a,\ v^i_{2a+2}=g-i-a+c-1 \ \ if\ \epsilon =1, \
                       v^i_{2c+2}=g-i-a+c \ \ if\ \epsilon =2$$
$$i=a^2+2a+1, u^i_{2a+1}=u^i_{2a+2}=i-1,\ v^i_{2a+1}=v^i_{2a+2}=g-i$$
$$u^i_{2e-1}=u^i_{2e}=i+e-2, \  v^i_{2e-1}=v^i_{2e}=g-i-e, e=a+2...k$$
We use Remark \ref{slce} to show that such a space of sections of
 the given vector bundle exists
(and is unique satisfying the conditions). We need to show that
the vanishing $u_{\alpha}, v_{\alpha}$ add up to $g-2$ except for
two values $\bar u_{\alpha}, \bar v_{\alpha}, \hat u_{\alpha},
\hat v_{\alpha}$ that add up to $g-1$ and then
$$E_i={ O}(\bar u_{\alpha}P_i+ \bar v_{\alpha}Q_i)\oplus
{ O}(\hat u_{\alpha}P_i+ \hat v_{\alpha}Q_i)$$

For $\epsilon =1 $,
$$u^i_{\alpha}+v^i_{\alpha}=g-2, \alpha\not= 2c+1, 2a+1\
u^i_{\alpha}+v^i_{\alpha}=g-1, \alpha = 2c+1, 2a+1$$
and
$$E_i={ O}(u^i_{2c+1}P_i+v^i_{2c+1}Q_i)\oplus
{ O}(u^i_{2a+1}P_i+v^i_{2a+1}Q_i)$$ Similarly, for $\epsilon =2$
$$u^i_{\alpha}+v^i_{\alpha}=g-2, \alpha\not= 2c+2, 2a+2\
u^i_{\alpha}+v^i_{\alpha}=g-1, \alpha = 2c+2, 2a+2$$
$$E_i={ O}(u^i_{2c+2}P_i+v^i_{2c+2}Q_i)\oplus
{ O}(u^i_{2a+2}P_i+v^i_{2a+2}Q_i)$$ Finally, for $c=a, \epsilon
=1$ (and hence $i=a^2+2a+1=(a+1)^2$)
$$u^i_{\alpha}+v^i_{\alpha}=g-2, \alpha\not= 2a+1,2a+2\
u^i_{\alpha}+v^i_{\alpha}=g-1, \alpha = 2a+1, 2a+2$$
$$E_i={ O}(u^i_{2a+1}P_i+v^i_{2a+1}Q_i)\oplus
{ O}(u^i_{2a+2}P_i+v^i_{2a+2}Q_i)$$ On the curves $C_i, i>k^2$,
take direct sum of an arbitrary line bundle of degree $g-1$ and
its complement with respect to the limit canonical series, namely
$$L_i\oplus L'_i,\ \ such \ that \ L_i\otimes L'_i=
{ O}(2(g-i)Q_i+2(i-1)P_i)$$
 There is a well determined subspace
of dimension $k$ of sections of $E_i$ with given vanishings $u^i_j$
at $P_i$ and $v^i_j=g-2-u^i_j$ at $Q_i$.
Take
$$u^i_{2e-1}=u^i_{2e}=i+e-k_1-2,
v^i_{2e-1}=v^i_{2e}=g-i+k_1-e,i=(k_1)^2+1...g-1,e=1...k_1$$

Take then the gluing at the nodes so that the section with
vanishing $v^{i-1}_j$ glues with the section with vanishing
$u^i_j$.

One can check that these spaces satisfy the conditions
 on the degree and vanishing of a limit linear series taking $a=g-1$ :
As every vector bundle is a direct sum of two line bundles of
degree $g-1=a$, the condition on the sums of degrees is satisfied and any section of
$E_i(-aP_i)$ or $E_i(-aQ_i)$ is determined by its value at the node.
Also, using the explicit expression for the $u^i_j,v^i_j$, one can check
that $v^{i-1}_j+u^i_j=a$ for all pairs $i,j$.

The given family of linear series is not part of a larger family
of linear series of dimension $k$ with fixed determinant $K$ : the
two ways in which one could make the vector bundle more general
are either by taking on some of the first $(k_1)^2$ components a
generic sum of two line bundles (with product equal to the
restriction of the canonical) or by taking more general gluings.
In either case, the vanishings $v^i_j$ at the corresponding
component $C_i$  will decrease. The values of the $v^{g-1}_j$ at
$C_{g-1}$ were
$$k,k,k-1,k-1,...2,2,1,1$$
If they were smaller, it would no longer be possible to glue this linear
series with two generic line bundles on $C_g$ and preserve the sections.

We now count the dimension of the family obtained in this way.

The gluing of $C_{a^2+2c}$ with $C_{a^2+2c+1}$ is free
because for every vanishing there is a two dimensional
space of sections with this vanishing
Therefore, these gluings
depend on four parameters. On the other hand, the gluing of
$C_{a^2+2c+1}$ with $C_{a^2+2c+2}$ must glue two given directions
to each other. Hence, it depends on only two parameters. The gluing between
$C_{a^2+2a}$ and $C_{a^2+2a+1}$ is again free.
So, between $C_{a^2}$ and $C_{(a+1)^2}$ there are $a+1$ gluings depending
on four parameters and $a$ gluings depending on two parameters.

The remaining $g-(k_1)^2$ gluings are free and depend on four parameters
each. The total contribution of the gluings is
$$4[(2+3+4+...+k_1)+g-(k_1)^2]+2(1+2+...+(k_1-1))=4g-(k_1)^2+k_1-4$$

The first $(k_1)^2$ vector bundles
are well determined while the remaining $g-(k_1)^2$ depend
on the choice of a line bundle of degree $g-1$ on the elliptic curve
(the second line bundle being its conjugate with respect to the
line bundle giving the canonical linear series).
 As the jacobian of an elliptic curve is one dimensional,
this contributes $g-(k_1)^2$ parameters in all.

We now need to compute the (negative) contribution of the endomorphisms.
The direct sum of two different line bundles has a two dimensional
family of endomorphisms.  The direct sum of two identical line
bundles has a four dimensional family of endomorphisms.
The vector bundles $E_1,E_4,E_9,...,E_{(k_1)^2}$ are of the latter type,
all of the remaining ones of the former type.

Hence, the dimension of the family is
$$4g-(k_1)^2+k_1-4+g-k_1^2-4k_1-2(g-k_1)+1$$
where the last +1 comes from the fact that the resulting vector bundle
is stable (see below) and therefore has a one dimensional family of endomorphisms.
The number above equals
$$3g-2(k_1)^2-k_1-3=3g-3-{k+1\choose 2}=\rho^k _K$$
as promised.

It remains to show that the vector bundle that we obtain in this way is
stable. As the restriction to each of the elliptic components is semistable,
using \cite{Tmfstc} or \cite{Tmfscr}, the resulting vector bundle is
 semistable.
We need to check that it is actually stable. It suffices then to show that
the line subbundles that make the restriction to the elliptic components
strictly semistable, don't glue with each other.
If $g\ge k_1^2+1$, there are precisely two subline bundles of $E_g$
that contradict stability. If $k_1=2,g=5$, there are two subline
bundles on the chain consisting of the first four elliptic curves
that contradict stability. But by the genericity of the gluing
between $C_4$ and $C_5$, these don't glue together.
If $k_1\ge 3$,
no line subbundles on the first $k_1^2$ components contradict
stability. This completes the proof of stability for the resulting
bundle.
\end{proof}
\bigskip

We now turn to the case of an odd number of sections $k=2k_1+1$
\begin{proof}
On the first $k_1$ elliptic components, one takes the same
vector bundles and spaces of sections with just
one more section that in the case $k=2k_1$. This new section
vanishes with multiplicities
$$u^i_{2k_1+1}=i+k_1-1,\  v^i_{2k_1+1}=g-i-k_1-1$$
so
$$u^i_{2k_1+1}+ v^i_{2k_1+1}=g-2$$
Hence
$$u^i_{2k_1+1}=u^i_{2k_1}+1, \ v^i_{2k_1+1}=v^i_{2k_1}-1,\ i=1...(k_1)^2$$
On the curves $C_i,\ i=k_1^2+m,\ 1\le m\le k_1$ we take the vector bundle to be
$${ O}((g-i-m+k_1+1)Q+(i+m-k_1-2)P)\oplus
{ O}((g-(k_1)^2-k_1-1)Q+((k_1)^2+k_1)P)$$ The vanishing $u^i_j,
v^i_j$ of the sections  of the linear series at the points
$P_i,Q_i$ in the curve $C_i,\ i=(k_1)^2+m,\ 1\le m\le k_1$ are
given by
$$u_{2e-1}=i+e-k_1-3, e< m,\ u_{2e-1}=i+e-k_1-2, m\le e\le k_1,
u_{2k_1+1}=(k_1)^2+k_1$$
$$u_{2e}=i+e-k_1-2, e\le k_1$$
$$v_{2e-1}=g-i-e+k_1+1, e\le m,\ v_{2e-1}=g-i-e+k_1, m<e\le k_1,$$
$$v_{2k_1+1}=g-1-(k_1)^2-k_,
v_{2e}=g-i-e+k_1, e\le k_1$$ On the curve $C_{(k_1)^2+k_1+1}$,
take the unique indecomposable vector bundle of rank two and
degree $2g-2$ which has one section vanishing with multiplicity
$g-1-(k_1)^2-k_1$ at Q and $(k_1)^2+k_1$ at P. On the remaining
curves, take again a direct sum of two line bundles each of degree
$g-1$ and whose product is the restriction of the canonical to
this component.

In order to compute the dimension of this family, note that on the first
$(k_1)^2$ components, the counts are the same as in the even number
of sections case. The gluing between the curve $C_{(k_1)^2}$ and
$C_{(k_1)^2+1}$ is free and therefore depends on four parameters.
The gluings between $C_{(k_1)^2+i}$ and $C_{(k_1)^2+i+1}$ must
glue together the direction corresponding to the vanishings
$v^{(k_1)^2+i}_{2k_1+1}$ and $u^{(k_1)^2+i+1}_{2k_1+1}$
Hence, they depend on three parameters.
All  vector bundles $E_{k_1^2+i},i=1...k_1+1$
 have a two dimensional family of endomorphisms.
Computing as before, we find that the dimension of the family is
$\rho^k _K$

Again, the vector bundle obtained in this form is stable
for $k_1\ge 2$ and for $k_1=1, g\ge 4$. For $k_1=1, g=3$
one can take as the vector bundle $E$ on any non-singular
non-hyperelliptic curve the dual of the kernel of the
evaluation map defined by
$$0\rightarrow E^*\rightarrow H^0(K)\otimes { O}_C
\rightarrow K\rightarrow 0$$
This is known to be stable \cite{PR}.

\end{proof}

\end{section}

\end{document}